\documentclass[12pt]{article}
\usepackage{a4wide}
\newcommand{\version}{October, 2011 }

\usepackage{amsthm,amsfonts,amsmath, amscd}
\swapnumbers
\theoremstyle{plain}

\newtheorem{thm}{THEOREM}[section]

\theoremstyle{definition}

\theoremstyle{definition}
\newtheorem{remark}[thm]{Remark}
\newcommand{\upchi}{\raise1pt\hbox{$\chi$}}
\newcommand{\R}{{\mathord{\mathbb R}}}

\newcommand{\cov}{{\rm cov}}
\newcommand{\var}{{\rm var}}
\newcommand{\dd}{{\, \rm d}}
\newcommand{\hess}{{\rm Hess}}
\renewcommand{\|}{{\Vert}}
\numberwithin{equation}{section}
\begin{document}

%\markboth{\scriptsize{CL \version}}{\scriptsize{CL \version}}

\def\tr{{\rm Tr}}

\title{Asymmetric Covariance Estimates of  \\Brascamp-Lieb  Type 
and \\ Related Inequalities for Log-concave Measures}
\author{\vspace{5pt} Eric A. Carlen$^1$, Dario
Cordero-Erausquin$^{2}$  and
Elliott H. Lieb$^{3}$ \\
\vspace{5pt}\small{$1.$ Department of Mathematics, Hill Center,}\\[-6pt]
\small{Rutgers University,
110 Frelinghuysen Road
Piscataway NJ 08854-8019 USA}\\
\vspace{5pt}\small{$2.$ Institut de Math\'ematiques de Jussieu, }\\[-6pt]
\small{Universit\'e Pierre et Marie Curie (Paris 6), 4 place Jussieu, 75252 Paris  France}\\
\vspace{5pt}\small{$3.$ Departments of Mathematics and
Physics, Jadwin
Hall,} \\[-6pt]
\small{Princeton University, P.~O.~Box 708, Princeton, NJ
  08542-0708}\\
 }
\date{\version}
\maketitle 
\footnotetext                                                                         
[1]{Work partially
supported by U.S. National Science Foundation
grant DMS 0901632.    }                           
                           \footnotetext
[2]{Work partially
supported by U.S. National Science Foundation
grant PHY 0965859.\\
\copyright\, 2011 by the authors. This paper may be
reproduced, in its
entirety, for non-commercial purposes.}

\begin{abstract}
An inequality of  Brascamp and Lieb  provides a
bound on the covariance of two functions with respect to
log-concave measures.  The bound estimates the covariance
by the product of the $L^2$ norms of the gradients of the
functions,
where the magnitude of the gradient is computed using an inner product
given by the inverse Hessian matrix of the potential of the log-concave measure.
Menz and Otto \cite{OM} proved a variant of this
with  the two $L^2$  norms replaced by $L^1$ and $L^\infty$
norms, but
only for $\R^1$. We prove a generalization of both by
extending these inequalities to $L^p$ and $L^q$ norms and
on $\R^n$, for any $n\geq 1$. We also prove an inequality
for integrals of divided differences of functions in terms
of integrals of their gradients.
\end{abstract}

\medskip
\leftline{\footnotesize{\qquad Mathematics subject
classification number:  26D10}}
\leftline{\footnotesize{\qquad Key Words: convexity,
log-concavity, Poincar\'e inequality}}

\section{Introduction} \label{intro}

 Let $f$ be a $C^2$ strictly convex function on $\R^n$ such
that $e^{-f}$ is integrable. By strictly convex, we mean
that the Hessian matrix, $\hess_f$, of $f$ is everywhere
positive.

Adding  a constant to $f$, we may suppose that
$$\int_{\R^n}e^{-f(x)}\dd^n x = 1\ .$$
Let $\dd \mu$ denote the probability measure
\begin{equation} \label{meas}
 \dd \mu := e^{-f(x)}\dd^n x \  ,
\end{equation}
 and let $\|\cdot\|_p$  denote the corresponding $L^p(\mu)$-norm.

For any two real-valued functions $f,g \in L^2(\mu)$, the covariance of $f$ and $g$ is the quantity
\begin{equation}\label{cov1}
\cov(g,h) := \int_{\R^n} gh\dd \mu - \left( \int_{\R^n} g\dd \mu  \right) \left( \int_{\R^n} h\dd \mu  \right)\ , 
\end{equation}
and the variance of $h$ is ${\rm var}(h) = \cov(h,h)$.

The Brascamp-Lieb (BL) inequality \cite{BL1}  for the variance of $h$ is
\begin{equation}\label{BL}
{\rm var}(h) \leq \int_{\R^n} ( \nabla h, \hess_f^{-1} \nabla h) \dd \mu\ ,
\end{equation}
where $(x,y)$ denotes the inner product in $\R^n$. (We shall also use $x\cdot y$ to denote this same inner product in simpler expressions
where it is more convenient.)

Since  $(\cov(g,h))^2\leq {\rm var }(g) {\rm var }(h)$, an
immediate consequence of (\ref{BL}) is
\begin{equation}\label{BL2}
(\cov(g,h))^2 \leq \int_{\R^n}  ( \nabla g, \hess_f^{-1}
\nabla g) \dd \mu
\int_{\R^n} ( \nabla h, \hess_f^{-1} \nabla h) \dd \mu\ .
\end{equation}

The one-dimensional variant of (\ref{BL2}), due to Otto and
Menz \cite{OM}, is
\begin{equation}\label{BL3}
|\cov(g,h)| \leq \| \nabla g\|_1
\| \hess_f^{-1} \nabla h\|_\infty^{\phantom{\int}}  \  
= 
\sup_x \left\{\frac{|h'(x)|}{f'' (x)}\right\}  
\  \int_\R | g'(x)| \dd \mu(x) 
\end{equation}
for functions $g$ and $h$ on $\R^1$. They  call this an {\em
asymmetric Brascamp-Lieb inequality}.
Note that it is asymmetric in {\it two} respects: One
respect is to  take an $L^1$
norm  of  $\nabla g$  and an $L^\infty$ norm of $\nabla h$,
instead of $L^2$ and $L^2$.  The second
respect is that the $L^\infty $ norm is weighted with the
inverse Hessian -- which here is
simply a number -- while the  $L^1 $ norm is not weighted.

Our first result is the following theorem, which generalizes
both (\ref{BL2}) and
(\ref{BL3}).

\begin{thm}[Assymetric BL inequality] 
\label{extott}
Let $\dd \mu (x) $ be as in \eqref{meas} and 
let $\lambda_{\rm min}(x)$ denote the least eigenvalue of
 $\hess_f(x)$.
For any locally Lipschitz functions $g$ and $h$ on $\R^n$ that are 
square integrable with respect to $\dd \mu$, and for
$2\leq p \leq \infty$, $1/p +1/q=1$, 
\begin{equation}\label{BL4}
|\cov(g,h)| \,\leq\, \big\Vert \hess^{-1/p}_f\, \nabla g \big\Vert_q \;
\big\Vert \lambda_{\rm min}^{(2-p)/p } \hess_f^{-1/p}\nabla h \big\Vert_p\ .
\end{equation}
This is sharp in the sense that (\ref{BL4}) cannot hold,
generally,  with a constant smaller than 1 on the
right side.
\end{thm}

%Note that for $n=1$, (\ref{BL4}) is (\ref{BL3})
%since for $n=1$, $\hess_f$ is a one by one matrix with entry
%$\lambda_{\rm min}$. 
For $p=2$, (\ref{BL4}) is (\ref{BL2}). 
Note that (\ref{BL4}) implies in particular that for Lipschitz functions $g,h$ on $\R^n$,
$$
|\cov(g,h)| \leq \big\Vert \lambda_{\rm min}^{-1/p}\, \nabla g \big\Vert_q\, \big\Vert \lambda_{\rm min}^{-1/q }\,  \nabla h \big\Vert_p\ .
$$
For $p=\infty$ and $q=1$, the latter is  
\begin{equation}\label{BL4.2}
|\cov(g,h)| \leq \big\Vert \nabla g\big\Vert_1\, \big\Vert \lambda_{\rm min}^{-1 } \nabla h \big\Vert_\infty^{\phantom{\int}} , 
\end{equation}
which for $n=1$ reproduces exactly (\ref{BL3}).
%For $p=2$, (\ref{BL4}) is (\ref{BL2}). 

We also prove the following theorem. In addition to its
intrinsic interest, it gives rise to an
alternative proof, which we give later, 
of Theorem~\ref{extott} in the case $p=\infty$ (though this proof 
only yields the sharp constant for $\R^1$, which is the
original
Otto-Menz case (\ref{BL3})).

\begin{thm}[Divided differences and gradients]\label{grbnd}
Let $\mu$ be  a probability measure with log-concave
density~\eqref{meas} . For any locally Lipschitz function $h$
 on $\R^n$,
\begin{equation}\label{t1}
 \int_{\R^n}\int_{\R^n} \frac {|h(x) - h(y)|}  {|x-y|} \dd
\mu(x)\dd \mu(y)  \leq 2^n\int_{\R^n}|\nabla h(x)|\dd \mu\ .
 \end{equation}
 \smallskip
\end{thm}

\iffalse
\begin{remark}The $m=1$ case is of particular interest.
Then, as we explain later, this gives us a bound on the
Cheeger constant for $\dd \mu$. To see this, Let $h$ be the
characteristic function of a set $A$. With this choice, and
with $m=1$, (\ref{t1B}) becomes

$$
2\mu(A) \left(1-\mu(A) \right) \leq C \int_{\partial A}
e^{-f} \dd \mathcal{H}_{n-1}  \   ,
$$
where $\mathcal{H}_{n-1}$ is $n-1$-dimensional Hausdorff
measure. Since the left side is not smaller than the minimum
of $\mu(A) $ and $\mu(A^c)$, the constant $C$ in inequality
(\ref{t1B}) is an upper bound for the Cheeger constant for
the measure
$\mu$. 
\end{remark}
\fi

\begin{remark} The constant $2^n$ is not optimal, as indicated by the examples in Section~\ref{examples} (we will actually briefly mention how to reach the constant $2^{n/2}$).
We do not know whether the correct constant grows
with $n$ (and then how), or is bounded uniformly in $n$. We do know that
for $n=1$, the constant is at least $2\ln 2$. We will return
to this later.
\end{remark}

The rest of the paper is organized as follows: Section 2 contains the proof of Theorem~\ref{extott}, and Section 3 contains the proof of
Theorem~\ref{grbnd}, as well as an explanation of the connection between the two theorems. Section 4 contains comments and examples
concerning the constant and optimizers in  Theorem~\ref{grbnd}. Section 5 contains a discussion of an application that motivated Otto and Menz,
and finally, Section 6 is an appendix providing some additional details on the original proof of the Brascamp-Lieb inequalities, which proceeds by induction on the dimension, and has an interesting connection with the application discussed in Section 5. 
\medskip

We end this introduction by expressing our gratitude to D. Bakry and M. Ledoux for fruitful exchanges on the preliminary version of our work. We originally proved~\eqref{BL4.2} with the constant $n2^n$ using Theorem~\ref{grbnd}, as explained in Section 3. Bakry and Ledoux pointed out to us that using a stochastic representation of the gradient along the semi-group associated to $\mu$ (sometimes referred to as the Bismut formula), one could derive inequality~\eqref{BL4.2} with the right constant $1$. This provided evidence that something more algebraic was at stake. It was confirmed by our general statement Theorem~\ref{extott} and by its proof below.

\section{Bounds on Covariance}

The starting point of the proof we now give for Theorem~\ref{extott}  is
a classical dual representation for the covariance which, in
the somewhat parallel setting of plurisubharmonic
potentials, goes back to the work of   H\"ormander. We shall
then adapt to our $L^p$ setting  H\"ormander's $L^2$
approach~\cite{hormander} to spectral estimates.

Let $g$ and $h$ be smooth and  compactly supported  on $\R^n$.  
Define the operator $L$ by
\begin{equation}\label{defL}
L = \Delta - \nabla f\cdot \nabla\ ,
\end{equation}
and note that
\begin{equation}\label{cobo1}
\int_{\R^n} g(x) L h(x) \dd \mu(x) = -\int_{\R^n}\nabla g(x)\cdot \nabla h(x) \dd \mu(x)\ ,
\end{equation}
so that $L$ is self-adjoint on $L^2(\mu)$.    Let us (temporarily)
add $\epsilon |x|^2$ to $f$ to make it uniformly convex, so that
the Hessian of $f$ is invertible and so that the operator $L$ has a
spectral gap. (Actually, $L$ always has a  spectral gap since $\mu$
is a log-concave probability measure, as noted in~\cite{KLS, Bo}.
Our simple regularization makes our proof independent of these deep
results.)

Then provided
\begin{equation}\label{cobo2}
\int_{\R^n} h(x) \dd \mu(x) = 0\ ,
\end{equation}
\begin{equation}\label{cobo4}
u := -\int_0^\infty e^{tL} h(x) \dd t
\end{equation}
exists and is in the domain of $L$, and satisfies $Lu =h$. 

Thus, assuming  (\ref{cobo2}), and by standard approximation arguments, 
\begin{eqnarray}\label{cobo3}
\cov(g,h) &=&   \int_{\R^n} g(x) h(x) \dd \mu(x)  = \int_{\R^n} g(x) L u(x) \dd \mu(x)\nonumber\\
&=& -\int_{\R^n}\nabla g(x)\cdot \nabla u(x) \dd \mu(x)\ .
\end{eqnarray}
This representation for the covariance is the starting point of the
proof we now give for Theorem~\ref{extott}.

\medskip
\noindent{\bf Proof of Theorem~\ref{extott}:}   Fix $2\leq p < \infty$, and let $q = p/(p-1)$, as in the statement of the theorem. 
Suppose $h$ satisfies (\ref{cobo2}), and define $u$ by (\ref{cobo4}) so that $Lu = h$. Then from (\ref{cobo3}),
\begin{eqnarray}\label{cobo5}
|\cov(g,h)| &\leq& \left| \int_{\R^n}\nabla g(x)\cdot \nabla u(x) \dd \mu(x) \right|\nonumber\\
&\leq&   \int_{\R^n} \big| \hess_f^{-1/p}\nabla g(x)\cdot \hess_f^{1/p}\nabla u(x) \big| \dd \mu(x)\nonumber\\
&\leq&  \|\hess_f^{-1/p}\nabla g(x)\|_q \, \| \hess_f^{1/p}\nabla u(x) \|_p\ .
\end{eqnarray}

Thus, to prove (\ref{BL4}) for $2\leq p < \infty$, it suffices to prove the following $W^{-1,p}$--$W^{1,p}$ type estimate: 
\begin{equation}\label{cobo7}
 \| \hess_f^{1/p}\nabla u(x) \|_p  \leq  \| \lambda_{\rm min}^{(2-p)/p } \hess_f^{-1/p}\nabla
h \|_p\ .
\end{equation}

Toward this end, we compute
\begin{eqnarray}\label{cobo8}
L(|\nabla u|^p) 
&=& p|\nabla u|^{p-2}(L\nabla u)\cdot \nabla u \nonumber \\
&\phantom{=}& \qquad +
p |\nabla u|^{p-2} \tr(\hess_u^2) + p(p-2)|\nabla
u|^{p-4}|\hess_u\nabla u|^2 
\nonumber\\
&\geq&  p|\nabla u|^{p-2}(L\nabla u)\cdot \nabla u\ ,
\end{eqnarray}
where we have used the fact that $p\geq 2$, and where the notation $L(\nabla u)$ refers to the  coordinate-wise action $(L\partial_1 u , \ldots ,L\partial_n u)$ of $L$.

Then, using the commutation formula (see the remark below)
\begin{equation}\label{comm}
 L(\nabla u) = \nabla(Lu) + \hess_f\nabla u\ ,
 \end{equation}
we obtain
$$0 = \int_{\R^n}L(|\nabla u|^p) \dd \mu(x) \geq  p\int_{\R^n}|\nabla u|^{p-2} \nabla u\cdot \nabla h\dd \mu(x)  +  
 p\int_{\R^n}|\nabla u|^{p-2} \nabla u\cdot  \hess_f \nabla u\dd \mu(x)\ , $$
 and hence
\begin{equation}\label{cobo9}
\int_{\R^n}|\nabla u|^{p-2} | \hess_f^{1/2} \nabla u|^2 \dd \mu(x)  \leq 
 \int_{\R^n}|\nabla u|^{p-2}| \hess_f^{1/p} \nabla u|| \hess_f^{-1/p} \nabla h|\dd \mu(x)\ .
 \end{equation}
We now observe that for any positive $n\times n$ matrix and any vector $v\in \R^n$,
$$ |A^{1/p}v|^p  \leq |v|^{p-2}|A^{1/2} v|^2 \ .$$
To see this, note that we may suppose $|v| =1$.  Then in the spectral representation of $A$, by Jensen's inequality,
$$|A^{1/p}v|  = \left(\sum_{j=1}^n \lambda_j^{1/p}v_j^2\right)^{1/2} \leq  \left(\sum_{j=1}^n \lambda_j^{1/2}v_j^2\right)^{1/p}\ .$$
Using this on the left side of (\ref{cobo9}), and using the obvious estimate
$$|\nabla u| \leq \lambda_{\rm min}^{-1/p}|\hess_f^{1/p}\nabla u|$$ 
on the right, we have
\begin{equation}\label{cobo10}
 \| \hess_f^{1/p} \nabla u\|_p^p \leq 
 \int_{\R^n}| \hess_f^{1/p} \nabla u|^{p-1} |\lambda_{\rm min}^{(2-p)/p} \hess_f^{-1/p} \nabla h|\dd \mu(x)\ .
 \end{equation}
 Then by H\"older's inequality we obtain (\ref{cobo7}).  
 
 It is now obvious that we can take the limit in which $\epsilon$ tends to zero, so that we obtain the inequality without any additional hypotheses on $f$.  Our calculations so far have required $2\leq p < \infty$, however, having obtained the inequality for such $p$, 
by taking the limit in which $p$ goes to infinity, we obtain the $p=\infty$, $q=1$  case of the theorem. 
 
 Finally, considering the case in which
 $$\dd \mu(x) =  (2\pi)^{-n/2}e^{-|x|^2}\dd x \ ,$$
 and $g = h = x_1$, we have that $\hess_f = {\rm I_d}$ and so
 $$\lambda_{\rm min} = |\hess_f^{-1/p} \nabla g| = |\hess_f^{-1/p}\nabla h| =1$$
 for all $x$, and so the constant is sharp, as claimed. \qed
 
 \begin{remark}  Many special cases and variants of the commutation relation (\ref{comm}) are well-known under different names. 
 Perhaps most directly relevant here is the case in which $f(x) = |x|^2/2$.  Then $\partial_j$ and its adjoint in $L^2(\mu)$,
 $\partial_j^* = x_j - \partial_j$, satisfy the canonical commutation relations, and the operator $L = -\sum_{j=1}^n \partial_j^*\partial_j$
 is (minus) the Harmonic oscillator Hamiltonian in the ground state representation. This special case of (\ref{comm}), in which the Hessian on the right is the identity, is the basis of the standard determination of the spectrum of the quantum harmonic oscillator using ``raising and lowering operators''. 
 
 In the setting of Riemannian manifolds, a commutation
relation analogous to (\ref{comm}) in which $L$ is the
Laplace-Beltrami operator and the Hessian is replaced by
${\rm Ric}$, the Ricci curvature tensor, is known as the
Bochner-Lichnerowicz formula. Both the Hessian version
(\ref{comm}) and 
 the Bochner-Lichnerowicz version have been used a number of
times to prove inequalities related to those we consider
here, for instance in 
 the work of Bakry and Emery on logarithmic Sobolev
inequalities.

We note that our proof immediately extends,
word for word, to the Riemannian setting if we use, in place
of~\eqref{comm} the  commutation satisfied by the operator
$L$ given by~(\ref{defL}) where $f$ is a (smooth) potential
on the manifold; That is, with some abuse of notation,
$L(\nabla u) = \nabla(Lu) + \hess_f\nabla u +
\textrm{Ric}\nabla u$, or rather, more rigorously, 
$$ 
L(|\nabla u|^p) \ge p|\nabla u|^{p-2} \big[ \nabla (Lu)
\cdot \nabla u + \hess_f \nabla u \cdot \nabla u +
\textrm{Ric} \nabla u \cdot \nabla u \, \big] .
$$
Thus, an analog of  Theorem~\ref{extott}
holds on a Riemannian manifold $M$ equipped with a
probability measure 
$$\dd\mu(x)=e^{-f(x)}\dd {\rm vol}(x)$$ 
where $ \dd {\rm vol}$ is the Riemannian element of volume
and $f$ a smooth function on $M$, provided $\hess_f$ at each
point $x$ is replaced in the statement by the
symmetric operator
$$H_x = \hess_f (x) + {\rm Ric}_x $$
defined on the tangent space.
Of course, the convexity condition on $f$ is accordingly
replaced by the  
assumption that $H_x>0$ at every point $x\in M$.
 
  \end{remark}

 \section{Bounds on Differences}
 
\noindent{\bf Proof of Theorem~\ref{grbnd}:}  Since $h(x) - h(y) = \int_0^1\nabla
h(x_t)\cdot (x-y)\dd t$, we have
\begin{equation}\label{gradbnd1}
 |h(x) - h(y)|  \leq  |x-y| \int_0^1|\nabla
h(x_t)|\dd t \quad{\rm where}\quad x_t := tx+(1-t)y\ .
\end{equation}
Next, by the convexity of $f$,
\begin{equation}\label{split}
e^{-f(x)}e^{-f(y)} =   e^{-(1-t)f(x)}e^{-tf(y)}
e^{-tf(x)}e^{-(1-t)f(y)}  \leq e^{-f(x_t)}
e^{-(1-t)f(x)}e^{-tf(y)} \ .
\end{equation}
Introduce the variables
\begin{eqnarray}\label{change}
w &=& tx + (1-t) y\nonumber\\
z &=& x-y\ .
 \end{eqnarray}
A simple computation of the Jacobian shows that this change
of variables is a measure preserving transformation for all
$0\leq t\leq 1$, and hence
\begin{multline}\label{divdi}
 \int_{\R^n}\int_{\R^n} \frac{|h(x) - h(y)|}{|x-y|}\dd
\mu(x)\dd \mu(y)  \leq \\  \int_0^1 \left(
\int_{\R^n}\int_{\R^n} |\nabla h(w)|
e^{-(1-t)f(w+(1-t)z)}e^{-tf(w-tz)}\dd z \, \dd
\mu(w)\right)\dd
t\ .
\end{multline} 
We estimate the right side of (\ref{divdi}).
By H\"older's inequality,
\begin{multline} \label{holder}
\int_{\R^n}e^{-(1-t)f(w+(1-t)z)}e^{-tf(w-tz)}\dd^n z \leq \\
\left(\int_{\R^n}e^{-f(w+(1-t)z)} \dd
z\right)^{1-t}\left(\int_{\R^n}e^{-f(w- tz)} \dd
z\right)^t \ .
\end{multline}
But
$$\int_{\R^n}e^{-f(w+(1-t)z)} \dd z = (1-t)^{-n}\qquad{\rm and}\qquad  \int_{\R^n}e^{-f(w- tz)} \dd z = t^{-n}\ ,$$
and finally,  $(1-t)^{-n(1-t)}t^{-nt}
= e^{-n(t\log t + (1-t)\log (1-t))}\leq 2^n
$. 
\hfill  \qed
\medskip

A corollary of Theorem \ref{grbnd} is a proof of Theorem
\ref{extott} for the special case of $q=1$ and $p=\infty$.
This proof is not only restricted to this case, it also has
the defect that the constant is not sharp, except in
one-dimension.  We give it, nevertheless, because it
establishes a link between the two theorems.

\medskip

\noindent{\bf Alternative Proof of Theorem~\ref{extott}
for $q=1$:} We shall use
the identity
\begin{equation}\label{cov2}
\cov(g,h)  =  \frac{1}{2}\int_{\R^n}\int_{\R^n}
[g(x)-g(y)][h(x)-h(y)] \dd \mu(x)\dd \mu(y)\ ,
\end{equation}
and estimate the differences on the right in different
ways.

Fix any $x\neq y$ in $\R^n$, and define the  vector $v
:=x-y$,
and for $0\leq t \leq 1$, define $x_t = y+tv  = tx+(1-t)y$.
Then for any Lipschitz function $h$,
\begin{equation}\label{cov3}
h(x) - h(y) = \int_0^{t} v\cdot \nabla h(x_t)\dd t\ .
\end{equation}

Now note that
\begin{equation}\label{pos1}
\frac{{\rm d}}{{\rm d}t} v\cdot \nabla f(x_t) = ( v,
\hess_f(x_t)v) \geq   |x-y|^2\lambda_{\rm min}(x_t) > 0\ .
\end{equation}
Integrating this in $t$ from $0$ to $1$, we obtain
\begin{equation}\label{pos2}
\left (x-y, \nabla f(x) - \nabla f(y) \right)= \int_0^1   
(v, \hess_f(x_t)v) \dd t >0\ ,
\end{equation}
which expresses the well-known monotonicity of gradients of
convex functions. 

Next, multiplying and dividing by $( v, \hess_f(x_t)v)$ in
(\ref{cov3}), we obtain
\begin{eqnarray}\label{cov67}
|h(x) - h(y)| &=&  \left|   \int_0^{1}  ( v, \hess_f(x_t)v)
( v, \hess_f(x_t)v)^{-1} v\cdot \nabla h(x_t)\dd
t\right|\nonumber\\
&\leq&    \int_0^{1}  ( v, \hess_f(x_t)v) \left|  ( v,
\hess_f(x_t)v)^{-1} v\cdot \nabla h(x_t)\right|\dd
t\nonumber\\
&\leq&    \int_0^{1}  ( v, \hess_f(x_t)v) \left|  \left(
\lambda_{\rm min}(x_t)\right)^{-1}|x-y|^{-2} v\cdot \nabla
h(x_t)\right|\dd t\nonumber\\
&\leq&   \sup_{z\in \R^n} \left\{   \frac{ |\nabla h(z)|}
{\lambda_{\rm min}(z)} \right\} |x-y|^{-1} \int_0^{1}  ( v,
\hess_f(x_t)v)\dd t\nonumber\\
&=&  \sup_{z\in \R^n} \left\{   \frac{ |\nabla h(z)|}
{\lambda_{\rm min}(z)} \right\} |x-y|^{-1}\left(x-y, \nabla
f(x) -  \nabla f(y)\right)\ .
\end{eqnarray} 

Define 
$$C :=  \sup_{z\in \R^n} \left\{   \frac{|\nabla h(z)|}
{\lambda_{\rm min}(z)} \right\} \ ,$$
and use (\ref{cov67}) in (\ref{cov2}):
\begin{eqnarray}
|\cov(g,h)| &\leq &   \frac{1}{2}\int_{\R^n}\int_{\R^n}
|g(x)-g(y) ||h(x)-h(y)| \dd \mu(x)\dd \mu(y)\nonumber\\
&\leq &   \frac{C}{2}\int_{\R^n}\int_{\R^n} |g(x)-g(y) |
\frac{1}{|x-y|}(x-y)\cdot \left[  \nabla f(x) -  \nabla
f(y)\right] \dd \mu(x)\dd \mu(y)\nonumber\\
&= &   \frac{C}{2}\int_{\R^n}\int_{\R^n} |g(x)-g(y) |
\frac{1}{|x-y|}(x-y)\cdot\left[  \nabla_y e^{-f(y)}e^{-f(x)}  - 
\nabla_x e^{- f(x)}e^{-f(y)} \right] \dd^n x\dd^n y\ .
\nonumber\\
&= &  -C\int_{\R^n}\int_{\R^n} |g(x)-g(y) |
\frac{1}{|x-y|}(x-y)\cdot  \nabla_x e^{- f(x)}e^{-f(y)} \dd^n
x\dd^n y\ ,\nonumber
\end{eqnarray}
where, in the last line,  we have used symmetry in $x$ and
$y$. 

Now integrate by parts in $x$.  Suppose first that $n>1$.
Then
$${\rm div}\left(\frac{1}{|z|}z\right) = \frac{n-1}{|z|}\
,$$
and
$|\nabla_x  |g(x)-g(y) | | = |\nabla_xg(x)|$ almost
everywhere.
Hence we obtain
\begin{equation}\label{cov6}
|\cov(g,h)| \leq C\left(\int_{\R^n}  |\nabla g(x)|\dd\mu(x)
+ (n-1) \int_{\R^n}\int_{\R^n} \frac {|g(x) - g(y)|} 
{|x-y|} \dd \mu(x)\dd \mu(y)\right)\ .
\end{equation}
For $n=1$,
${\displaystyle {\rm div}\left(\frac{1}{|z|}z\right) =
2\delta_0(z)}$
and (\ref{cov6}) is still valid since $|g(x) -
g(y)|\delta_0(x-y) = 0$.

Now, for $n=1$, (\ref{cov6}) reduces directly  to
(\ref{BL3}). For $n>1$, it reduces to~(\ref{BL4.2}) upon application of Theorem 1.2, but with the constant $n2^n$ instead of $1$.  \qed

\section {Examples and Remarks on Optimizers in Theorem
\ref{grbnd} \label{examples}}

Our first examples address the question of the importance of
log-concavity.

{\bf (1.) Some restriction on $\mu$ is necessary:} If a
measure $d\mu(x) =F(x) dx$ on $\R$ has $F(a)=0 $ for some $a\in \R$, and
$F$ has positive mass to the left and right of $a$, then
inequality \eqref{t1} cannot possibly hold with any
constant. The choice of $h$ to be the Heaviside step
function shows that \eqref{t1} cannot hold with any constant
for this $\mu$.

\medskip

{\bf (2.)  Unimodality is not enough:} Take  $d\mu(x) =F(x)
dx$, with $F(x) =
1/4\varepsilon$ on $(-\varepsilon, \varepsilon)$ and $F(x) =
1/4(1-\varepsilon )$ otherwise on the
interval $(-1,1)$ and $F(x) =0$ for $|x| >1$.  Let $g(x) =1$
for $|x| < \varepsilon +\delta$ and $g(x) =0$ otherwise.
When $\delta$ is positive but small, 
$$\int_{\R}|\nabla g|\dd \mu(x) = 1/2(1-\varepsilon )$$
while 
$$\int_{\R}\int_{\R}\frac{|g(x)-g(y)|}{|x-y|}\dd \mu(x)\dd
\mu(y) = O(-\ln(\epsilon))\ .$$

\medskip
{\bf (3.)  For $n=1$, the best constant in (\ref{t1}) is at
least $2\ln2$:}
Take  $d\mu(x) =F(x) dx$, with $F(x) = 1/2$ on $(-1,1)$  and
$F(x) =0$ for $|x| >1$.
Let $g(x) =1$ for $x\geq 0$ and $ g(x) = 0$ for $x< 0$. All
integrals are easily computed.

\medskip

{\bf (4.)  The best constant is achieved for characteristic
functions:} When seeking the best constant 
in (\ref{t1}),
it suffices, by a standard truncation argument,
 to consider bounded Lipschitz functions $h$. Then, since
neither side of the inequality is affected if we add a
constant to $h$, it suffices to consider non-negative
Lipschitz functions.    We use the layer-cake representation
\cite{LL}: 
 $$h(x) = \int_0^\infty \chi_{\{ h> t\} }(x) \dd t\ .$$
Then 
\begin{equation}\label{layer}\int_{\R^n}\int_{\R^n}\frac{
|h(x)-h(y)|}{|x-y|}\dd \mu(x)\dd \mu(y) \leq 
\int_0^\infty
\int_{\R^n}\int_{\R^n}\frac{|\chi_{\{h>t\}}(x)-\chi_{\{h>t\}
}(y)|}{|x-y|}
\dd \mu(x)\dd \mu(y)\dd t\end{equation}
Define $C_n$ to be the best constant for characteristic
functions of sets $A$ and log-concave measures $\mu$:
\begin{equation} \label{see}
C_n := \sup_{f, A} \left\{
\frac{\int_{\R^n}\int_{\R^n}\frac{|\chi_A(x)-\chi_A(y)|}{
|x-y|}\dd \mu(x)\dd \mu(y)}{\int_{\partial A}e^{-f(x)}\dd
{\mathcal H}_{n-1}(x)}\right\}
\end{equation}
where $ {\mathcal H}_{n-1}$ denotes $n-1$ dimensional
Hausdorff measure.
Apply this to (\ref{layer}) to conclude that
\begin{eqnarray} \label{cheeger}
\label{layer2}
\int_{\R^n}\int_{\R^n}\frac{|h(x)-h(y)|}{|x-y|}\dd \mu(x)\dd
\mu(y) &\leq&  C_n
\int_0^\infty\int_{\partial \chi_{\{h>t\}}}e^{-f(x)}\dd
{\mathcal H}_{n-1}(x)\dd t\nonumber\\
&=& C_n \int_{\R^n}|\nabla h(x)|\dd \mu(x)\ ,
\end{eqnarray} 
where the co-area formula was used in the last line.  Thus,
inequality
(\ref{t1}) holds with the constant $C_n$; in short, it
suffices to
consider characteristic functions as trial functions. 
Note that the argument is also valid at the level of each measure $\mu$ individually, although we are interested here in uniform bounds.

With characteristic functions in mind, let us consider the
case that
$g$ is the characteristic function of a half-space in
$\R^n$. Without
loss of generality let us take this to be $\{x\,:\, x_1
<0\}$.
Clearly, the left side of \eqref{t1} is less than the
integral with
$|x-y|^{-1} $ replaced by $|x_1-y_1|^{-1}$. Since the
marginal
(obtained by integrating over $x_2, \dots, x_n$) of a log
concave
function is log concave, we see that our inequality reduces
to the
one-dimensional case. In other words, the constant $C_n$  in
\eqref{see} would equal  $C_1$,
independent of $n$, if the supremum were restricted to
half-spaces instead of to arbitrary
measurable sets.

\medskip

{\bf (5.) Improved constants and geometry of log-concave
measures:}
With additional assumptions on the measure one can see that
the constant is not only bounded in $n$, but of order
$1/\sqrt n$.
We are grateful to F.~Barthe and M.~Ledoux for discussions and
improvements in particular cases 
concerning the constant in Theorem \ref{grbnd}. This relies on
the Cheeger constant $\alpha(\mu)^{-1}>0$ associated to the log-concave probability measure $\dd \mu$, which is defined to be the best constant in the inequality 
$$\forall A\subset \R^n \textrm{(regular enough)}, \qquad
\mu(A) (1- \mu (A))
\leq \alpha(\mu) 
\int_{\partial A}e^{-f(x)}\dd
{\mathcal H}_{n-1}(x)
$$
M. Ledoux suggested the following procedure. Split the 
function $|x-y|^{-1} $ into two pieces according to
whether $|x-y|$ is less than or greater than $R$, for some
$R>0$. With  $h$ being the characteristic function of $A$, the
contribution to the left side of \eqref{cheeger} for
$|x-y|>R$ is bounded above by $2R^{-1}  \alpha(\mu)
\int_{\partial A}e^{-f(x)}\dd
{\mathcal H}_{n-1}(x)$. The contribution for $|x-y| \leq R$
is bounded above in the same manner as in the proof of 
Theorem \ref{grbnd}, but this time we only have to integrate
$z$ over the domain $|z|\leq R$ in each of the integrals in
\eqref{holder}. Thus, our bound $2^n$ is improved by a 
factor, which is the  $\dd \mu$ volume of the ball
$B_R=\{|z|\leq R\}$, once we used the Brunn-Minkowski inequality for the bound
$$\mu\big((1-t) B_R+w)\big)^{1-t} \, \mu\big(t B_R + w\big)^t \le \mu\big(B_R + w\big) \le \overline{\mu}(B_R):=\sup_{x}\mu(B_R+x).$$
The final step is to optimize the sum of
the contributions of the two terms with respect to $R$. 
Thus, if we denote $C_n(\mu)$ the best constant in the inequality~\eqref{t1} of Theorem \ref{grbnd} for a fixed measure $\mu$, we have
\begin{equation}\label{Lbound}
C_n(\mu) \le \inf_{R>0}   \big\{ 2^n \overline{\mu} (B_R) + 2R^{-1} \alpha(\mu) \big\} \le 2^n.
\end{equation}
Note that if $\mu$ is symmetric (i.e. if $f$ is even), then the Brunn-Minkowski inequality ensures that $\overline{\mu}(B_R)=\mu(B_R)$.

Unlike in \eqref{t1}, this improved bound depends on $\mu$ but there are situation where this gives optimal estimates as pointed out to us by F.~Barthe. As an example, consider the case where $\mu$ is the standard Gaussian measure on $\R^n$. Using the known
value of the Cheeger constant for this $\mu$, and linear
trial functions, one finds that the constant is bounded
above and below by a constant times $n^{-1/2}$.
% for the Gaussian measure.

%\medskip

Actually, we can use (\ref{Lbound}) to 
improve the constant from $2^n$ to $2^{n/2}$ for arbitrary
measures using some recent results from the geometry
of log-concave measures. Without loss of generality, we can assume, by translation of $\mu$, that
 $ \int |x|\, d\mu(x)=\inf_v \int |x+v|\, d\mu(x)=:M_\mu $.  It was proved
in~\cite{KLS, Bo} that for every  log-concave
measure on $\R^n$,
$$\alpha(\mu) \le c M_\mu$$  
where $c>0$ is some numerical constant (meaning a possibly
large, but computable, constant, in particular independent
of $n$ and $\mu$, of course). On the other hand, it was
proved by Gu\'edon~\cite{Gue} that for every  log-concave
measure $\nu$ on $\R^n$
$$\nu(B_R) \le \frac{C}{\int |x|\, d\nu}\, R$$
for some numerical constant $C>0$. In the case $\mu$ is not symmetric, we pick  $v$  such that $\mu(B_r + v)=\overline{\mu}(B_R)$, and then we apply the previous bound to $\nu(\cdot)=\mu(\cdot +v)$ in order to get that $\overline{\mu}(B_R) \le \frac{C}{M_\mu}$. Using these two estimates in~\eqref{Lbound} we see that
$$C_n(\mu) \le \inf_{s>0} \{C 2^n s + c/s\} = \kappa \,  2^{n/2}$$
for some numerical constant $\kappa >0$.

The Brascamp-Lieb
inequality~\eqref{BL}, as well as
inequality~\eqref{t1}, have connections with the geometry of
convex bodies. It was observed in~\cite{BL2} that
\eqref{BL} can be deduced from the Pr\'ekopa-Leindler
inequality (which is a functional form of the
Brunn-Minkowski inequality). But the converse is also true:
the Pr\'ekopa theorem follows, by a local computation, from
the Brascamp-Lieb inequality (see~\cite{CE} where the
procedure is explained in the more general complex setting).
To sum up, the Brascamp-Lieb inequality~\eqref{BL} can be
seen as the \emph{local} form of the Brunn-Minkowski
inequality for convex bodies.

\section{Application to Conditional Expectations}

Otto and Menz were motivated  to prove
(\ref{BL3}) for an application that involves a large amount of additional structure that we cannot go into here. We shall however give an application of Theorem ~\ref{extott}  to a type of estimate
that is related to one of the central estimates in \cite{OM}.

We use the notation in \cite{BL1}, which is adapted to
working with a partitioned set of variables.  Write a point
$x\in \R^{n+m}$ as $x=(y,z)$ with $y\in \R^m$ and $z\in
\R^{n}$. For a function $A$ on $\R^{n+m}$, let $\langle
A\rangle_z(y)$
denote the conditional expectation of $A$ given $y$, with
respect to
$\mu$. For a function $B$ of $y$ alone, $\langle B\rangle_y$
is the
expected value of $B$, with respect to $\mu$.  As in
\cite{BL1}, a
subscript $y$ or $z$ on a function denotes differentiation
with respect to $y$ or $z$, while a subscript $y$ or $z$ on
a bracket denotes integration.
For instance, for a function $g$ on $\R^{n+m}$, $g_y$
denotes  the vector $\big(\frac{\partial g}{\partial
{y_i}}\big)_{i\le n}$ in $\R^n$, and for $i\le n$,
$g_{y_i z}$ denotes 
the vector $\big(\frac{\partial^2 g}{\partial
{y_i}\partial{z_j}}\big)_{j\le m}$  in $\R^m$. Finally, 
$(g_{yz})$ denotes the $n\times m$ matrix having the
previous vectors as rows.

Let $h$ be non-negative with $\langle h\rangle_x = 1$ so
that $h(x)\dd \mu(x)$
is a probability measure, and so is $\langle
h\rangle_z(y)\dd\nu(y)$, where $\dd \nu(y)$
is the marginal distribution of $y$ under $\dd\mu(x)$.  

A problem that frequently arises \cite{BH,GOVW,LPY,OR,OM} is
to  estimate the Fisher information of $\langle
h\rangle_z(y)\dd\nu(y)$ in terms of the
Fisher information of  $h(x)\dd \mu(x)$ by proving an
estimate of the form
\begin{equation}\label{appl}
\left\langle \frac{|(\langle h\rangle_z)_y|^2}{\langle
h\rangle_z}\right \rangle_y \leq 
C\left\langle \frac{|h_x|^2}{h}\right\rangle_x\ . 
\end{equation}

Direct differentiation under the integral sign in the variable $y_i$  gives
$$(\langle h\rangle_z)_{y_i} = \langle h_{y_i}\rangle_z - \cov_z(h,f_{y_i})\ ,$$
where $\cov_z$ denotes the conditional covariance of
$h(y,z)$
and $f_{y_i}(y,z)$, integrating in $z$ for each fixed $y$. 
Let $u = (u_1,\dots,u_m)$ be any unit vector in $\R^m$.  Then
Hence, for each $y$,
\begin{eqnarray}
(\langle h\rangle_z)_y\cdot u = \sum_{i=1}^m  (\langle h\rangle_z)_{y_i}u_i &= &
 \sum_{i=1}^m \langle h_{y_i}\rangle_z u_i -  \sum_{i=1}^m
 \cov_z(h,f_{y_i})u_i\nonumber\\
 &=&  \langle h_{y}\rangle_z\cdot u - \cov_z(h,f_{y}\cdot u)\ ,\nonumber
 \end{eqnarray}
 and hence, choosing $u$ to maximize the left hand side, 
\begin{equation}\label{eig2}
|(\langle h\rangle_z)_{y}|^2
\leq 2|\langle h_{y}\rangle_z|^2 +
2\left(\cov_z(h,f_{y}\cdot u)\right)^2\ .
\end{equation}

By (\ref{BL4}),
\begin{equation}\label{eig3}
|\cov_z(h,f_{y}\cdot u)| \leq
\langle|h_z|\rangle_z
\|\lambda_{\rm min}^{-1}|(f_{y}\cdot u)_z|\|_\infty^{\phantom{int}}\ . 
\end{equation}

Note that the  least eigenvalue of the $n\times n$
block $f_{zz}$ 
is at least as large as the least eigenvalue
$\lambda_{\rm min}(y,z)$ of the full 
Hessian, by the variational principle.  
Hence, while we are entitled to use
the least eigenvalue of the
$n\times n$ block $f_{zz}$ of the full $(n+m)\times(n+m)$
Hessian matrix $f_{xx}$, and this would be important in the application
in the one dimensional case made in \cite{OM}, here, without any special structure to take advantage of,  we simply use the least eigenvalue of the full matrix in our bound. 

Next note that
$$|(f_{y}\cdot u)_z|^2  \leq \sum_{i=1}^m\left(\sum_{j=1}^n (f_{y_i,z_j})^2\right)u_i^2\ ,$$
and that ${\displaystyle \sum_{j=1}^n (f_{y_i,z_j})^2}$
is the
$i,i$ entry of $f_{yz}^Tf_{yz}$ where $f_{yz}$ denotes the
upper right corner block of the Hessian matrix. This number
is no greater than the $i,i$ entry of the square of the
full Hessian matrix. This, in turn, is no greater than
$\lambda_{\rm max}^2$.  Then, since $u$ is a unit vector, we have
$$|(f_{y}\cdot u)_z| \leq  \lambda_{\rm max}\ .$$
Using this in (\ref{eig3}), we obtain
\begin{equation}\label{eig37}
|\cov_z(h,f_{y}\cdot u)| \leq
\langle|h_z|\rangle_z
\|\lambda_{\rm max}/\lambda_{\rm min}\|_\infty^{\phantom{int}}\ ,
\end{equation}
and then from  (\ref{eig2})
\begin{equation}\label{eig4}|(\langle h\rangle_z)_y|^2
\leq 2|\langle h_y\rangle_z|^2 +
2
\|\lambda_{\rm max}/\lambda_{\rm min}\|_\infty^{2}
\langle|h_z|\rangle_z^2\ .
\end{equation}
Then the Cauchy-Schwarz inequality yields
\begin{equation}\label{eig5}
\left(\langle|h_z|\rangle_z\right)^2 \leq 
\left\langle
\frac{|h_z|^2}{h}\right\rangle_z\langle h\rangle_z\ .
\end{equation}
Use this in (\ref{eig4}), divide both sides by $\langle
h\rangle_z$, and integrate in $y$. The joint convexity in
$A$ and $\alpha >0$ of $A^2/\alpha$ yields (\ref{appl})
with the constant $C = 2\|\lambda_{\rm max}/\lambda_{\rm min}\|_\infty^2$.

The bound we have obtained becomes useful when
$\lambda_{\rm max}(x)/\lambda_{\rm min}(x)$ is bounded
uniformly. 
Suppose that $f(x)$ has the form $f(x) =
\varphi(|x|^2)$. Then the eigenvalues of the Hessian of $f$
are $2\varphi'(|x|^2)$, with multiplicity $m+n-1$, and
$4\varphi''(|x|^2)|x|^2 + 2\varphi'(|x|^2)$, with
multiplicity 1. Then both eigenvalues are positive, and the
ratio is bounded, whenever $\varphi'$ is positive and, for
some $ c < 1 < C < \infty$,
$$ -c\varphi'(s) \leq s\varphi''(s) \leq
C\varphi'(s)\ .$$

\medskip

\begin{remark}[Other asymmetric variants of the BL inequality]
Together, (\ref{eig3}) and
(\ref{eig5}) yield
$$\frac{\left(\cov_z(h,f_y\cdot u)\right)^2}{ \langle h\rangle_z} 
\leq
\left\langle
\frac{|h_z|^2}{h}\right\rangle_z
\|\lambda_{\rm min}^{-1}(f_{y}\cdot u)_z)|\|_\infty^2\  .
$$
A weaker inequality is 
\begin{equation}\label{eig8}
\frac{\left(\cov_z(h,f_y\cdot u)\right)^2}{ \langle h\rangle_z}
\leq
\left\langle
\frac{|h_z|^2}{h}\right\rangle_z
\|\lambda_{\rm min}^{-1}\|_\infty^2 \|(f_{y}\cdot u)_z)\|_\infty^2\  .
\end{equation}

In the context of the application in \cite{OM},
finiteness of   $\|(f_{y}\cdot u)_z\|_\infty^{\phantom{int}}$
limits  $f$ to quadratic growth at infinity. A major
contribution of \cite{OM} is to remove this limitation in
applications of (\ref{appl}). 
The success of this application of (\ref{BL3}) depended on
the full weight of the inverse Hessian being allocated to
the $L^\infty$ term.

Nonetheless, once the topic of asymmetric BL inequalities is
raised, one might enquire whether an inequality of the type
\begin{equation}\label{BL35}
|\cov(g,h)| \leq C\| \nabla g\|_\infty^{\phantom{int}}
\| \hess_f^{-1} \nabla h\|_1\ 
\end{equation}
can hold for any constant $C$. 
There is no such inequality, even in one dimension. To see
this, suppose that for some $a\in \R$ and some
$\epsilon>0$, $f_{xx}> M$ on $(a-\epsilon,a+\epsilon)$. 
Take $h(x) = 1$ for $x>a$ and $h(x) = 0$ for $x\leq a$. 
Take $g(x) = x-a$.  Suppose that $f$ is even about $a$. Then
$\cov(g,h)   = \int_a^\infty (x-a)e^{-f(x)}\dd x$, while
$\| \hess_f^{-1} \nabla h\|_1 \leq M^{-1}$, and $f$ can be
chosen to make $M$ arbitrarily large while keeping
$\| \nabla g\|_\infty \leq 1$, and  $\cov(g,h) $ bounded
away from zero.

\end{remark}

\section{Appendix} We recall that the original proof of (\ref{BL}),
Theorem 4.1 of \cite{BL1}, used dimensional induction, though interesting non-inductive proofs have since been provided \cite{BL2}. 

The starting point for the inductive proof is that the proof for
$n=1$ is elementary. The proof of the inductive step is more involved,
and we take this opportunity to provide more detail about the passage
from eq.~(4.9) of \cite{BL1} to eq.~(4.10) of \cite{BL1}. There is an
interesting connection with the application discussed in
the previous section, which also concerns
$\langle h_y\rangle_z -
\cov_z(h,f_y)$. We continue using the notation introduced there, but now $m=1$ (i.e. $y\in \R$). 

Eq. (4.9) reads $ \var(h) \leq \langle B\rangle_y$ where
\begin{equation}\label{ind0}
B = {\rm var}_z(h) + \frac{ [\langle h_y\rangle_z -
\cov_z(h,f_y)]^2}
{\langle  f_{yy}\rangle_z - \var_z f_y}\ .
\end{equation}
Our goal is to prove 
\begin{equation} \label{goal}
B \leq \langle (h_z,f_{zz}^{-1}h_z)\rangle_z + \frac{ \langle h_y - 
(h_z,f_{zz}^{-1}f_{yz})\rangle_z^2}
{\langle  f_{yy} -(f_{yz},f_{zz}^{-1}f_{yz})\rangle_z}\ .
\end{equation}

To do this,  use the inductive hypothesis; i.e., for any $H$ on $\R^{n-1}$,
\begin{equation}\label{ind}
 \var_z (H) \leq \langle H_z,f_{zz}^{-1}H_z\rangle_z\ .
\end{equation}

Apply this to arbitrary linear combination $H = \lambda h + \mu f_y$ to
conclude the $2\times 2$ matrix inequality
$$
\left[\begin{array}{cc} \var_z (h) & \cov_z(h,f_y)\\ \cov_z(h,f_y) & \var_z(f_y)\end{array}\right] \leq
\left[\begin{array}{cc} \langle(h_z,f_{zz}^{-1}h_z)\rangle_z & \langle(h_z,f_{zz}^{-1}f_{yz})\rangle_z\\ 
\langle(f_{yz},f_{zz}^{-1}h_z)\rangle_z & \langle(f_{yz},f_{zz}^{-1}f_{yz})\rangle_z\end{array}\right]
$$
Take the determinant of the difference to find that
\begin{equation}\label{ind2}
\langle (h_z,f_{zz}^{-1}h_z)\rangle_z - \var_z(h) \geq \frac{[\langle (h_z,f_{zz}^{-1}f_{yz})\rangle_z - \cov_z(h,f_y)]^2}
 {\langle(f_{yz},f_{zz}^{-1}f_{yz})\rangle_z - \var_z(f_y)}\ .
\end{equation}
Combine \eqref{ind0} and \eqref{ind2} to obtain
\begin{equation}\label{ind3}
B \leq \langle (h_z,f_{zz}^{-1}h_z)\rangle_z  +  \frac{ [\langle h_y\rangle_z - \cov_z(h,f_y)]^2}
{\langle  f_{yy}\rangle_z - \var_z( f_y)}  - 
\frac{[\langle (h_z,f_{zz}^{-1}f_{yz})\rangle_z - \cov_z(h,f_y)]^2}
 {\langle(f_{yz},f_{zz}^{-1}f_{yz})\rangle_z - \var_z(f_y)}
\end{equation}
Since $a^2/\alpha$ is jointly convex in $a$ and $\alpha > 0$, and is homogeneous of degree one,
for all $\alpha > \beta > 0$ and all $a$ and $b$, 
$$\frac{a^2}{\alpha} \leq \frac{b^2}{\beta} + \frac{(a-b)^2}{\alpha-\beta}\ .$$
That is, $a^2/\alpha - b^2/\beta \leq
(a-b)^2/(\alpha-\beta)$. Use this on the right side
of \eqref{ind3} to obtain (\ref{goal}), noting that the
positivity of
$\alpha - \beta = 
\langle  f_{yy}\rangle_z - 
\langle(f_{yz},f_{zz}^{-1}f_{yz})\rangle_z$ is
a consequence of the positivity of the Hessian of $f$.

\end{document}